\documentclass{article}
\usepackage{amssymb}
\usepackage{amsmath}
\begin{document}

\newtheorem{theorem}{Theorem}
\newtheorem{proposition}[theorem]{Proposition}
\newtheorem{corollary}[theorem]{Corollary}
\newtheorem{lemma}[theorem]{Lemma}
\newtheorem{definition}[theorem]{Definition}
\newcommand{\QED}{\phantom{m}\hfill $\Box$}
\newcommand{\real}{\ensuremath{\Bbb{R}}}
\newcommand{\realn}{\ensuremath{\Bbb{R}^n}}
\newcommand{\realm}{\ensuremath{\Bbb{R}^m}}
\newcommand{\realk}{\ensuremath{\Bbb{R}^k}}
\newcommand{\nat}{\ensuremath{\Bbb}{N}}
\newcommand{\fet}[1]{\ensuremath{\mathbf{#1}}}

\title{\vspace*{-1.5cm}
\begin{flushright}
\begin{minipage}{4cm}
\tiny
{\sc Dept. of Math./CMA \hfill University of Oslo\\
Pure Mathematics \hfill     No 11\\
ISSN 0806--2439 \hfill     August 2007}
\end{minipage}
\end{flushright}
\vskip1cm
A weighted random walk approximation to fractional Brownian motion}
\author{Tom Lindstr\o m\thanks{Centre of Mathematics for Applications 
and Department of Mathematics, PO Box 1053 Blindern, N-0316 Oslo, 
Norway. e-mail:lindstro@math.uio.no}}
\date{}
\maketitle
\begin{abstract}
\noindent We present a random walk approximation to fractional
Brownian motion where the increments of the fractional random
walk are defined as a weighted sum of the past increments of
a Bernoulli random walk.\\



\end{abstract}

\noindent The purpose of this brief note is to describe a discrete approximation to fractional
Brownian motion. The approximation works for all Hurst indices $H$, but take slightly different
forms for $H\leq\frac{1}{2}$ and $H>\frac{1}{2}$. There are already several discrete
approximations to fractional Brownian motion in the literature (see, e.g., \cite{Taqqu}, \cite{Joseph}, \cite{Dasgupta}, \cite{Sottinen}, \cite{KK}, \cite{BCF}, \cite{KS}, \cite{Neuenkirch} for this and related topics), and the advantage of 
the present approach is that the increments of the fractional random walk is given as a weighted
sum of past increments of an ordinary (Bernoulli) random walk. This gives an excellent understanding 
of the dynamics of the process and is a good starting point for stochastic calculus with
respect to fractional Brownian motion. A similar idea is exploited in much greater generality by Konstantopoulos and Sakhanenko in  \cite{KS}, but they assume that $H>\frac{1}{2}$, while the present paper is mainly of interest when $H<\frac{1}{2}$. 

The discrete approximation is based on Mandelbrot and Van Ness' \cite{MVN} moving frame representation of fractional Brownian motion:

\[x_{t}=c_{H}\int_{-\infty}^t\left((t-r)^{H-\frac{1}{2}}-(-r)_{+}^{H-\frac{1}{2}}\right)\;db_{r}\]
where the scaling constant $c_{H}$ is given by
\[c_{H}=\left(\int_{0}^{\infty}\left((1+u)^{H-\frac{1}{2}}-u^{H-\frac{1}{2}}\right)^2\;du
+\frac{1}{2H}\right)^{-\frac{1}{2}}=\frac{\sqrt{\Gamma(2H+1)\sin(\pi H)}}{\Gamma(H+\frac{1}{2})}\]
(see also \cite{ST}). This representation will be used to establish the convergence.

\section{The main theorem}

To state the main result, we need some notation. For each natural number $N$, let 
$\Delta t_{N}=\frac{1}{N}$ and think of
\[T_{N}=\{k\Delta t_{N}\;|\;k\in\mathbb{Z}\}\]
as a timeline. We let $T_{N}^+$ denote the nonnegative part of $T$.
 It is convenient to use the following convention for sums over elements in $T_{N}$:
\[\sum_{r=s}^t f(r)=f(s)+f(s+\Delta t)+\cdots +f(t-\Delta t_N)\]
Note that the lower limit $s$ is included in the sum, but the upper limit $t$ is not. 
We shall also write $\Delta f(t)=f(t+\Delta t_{N})-f(t)$ for the forward increment of $f$ at $t$. 

For all $t\in T_{N}$, let $\omega_{N}(t)$ be independent random variables taking values
 $\pm 1$ with probability $\frac{1}{2}$. We shall write 
 $\Delta B_{N}(t)=\sqrt{\Delta t_{N}}\omega_{N}(t)$ and think of $B_{N}$ as a Bernoulli
 random walk approximating Brownian motion. For $0<H<1$ and $N\in\nat$, define a 
 process $X_{H,N}:\Omega_N\times T_{N}^+\to\real$ by $X_{N,H}(0)=0$ and 
 \[\Delta X_{H,N}(s)=K_{H}\Delta t_N^{H-\frac{1}{2}}\Delta B_N(s)+
\sum_{r=-\infty}^s (H-\frac{1}{2})(s-r)^{H-\frac{3}{2}}\Delta t_N\Delta B_N(r)\]
(using, e.g., Kolmogorov's one series theorem, see \cite{Var}, one easily checks that the sum converges a.s.) where the constant $K_{H}$ is defined by
\[K_{H}=\left\{\begin{array}{cl}-(H-\frac{1}{2})\zeta(\frac{3}{2}-H)&\mbox{ for }H<\frac{1}{2}\\
&\\
1&\mbox{ for }H\geq\frac{1}{2}\end{array}\right.\]
(as usual, $\zeta(s)=\sum_{n=1}^{\infty}n^{-s}$ when $s>1$). Except for the Mandelbrot-Van Ness scaling factor $c_H$, 
$X_{H,N}$ will be our random walk approximation to fractional Brownian motion. For convergence puposes it will be convenient to think of $X_{H,N}$ as a c\`{a}dl\`{a}g process defined on $[0,\infty)$, and we do this simply by assuming that $X_{H,N}$ is constant between points in $T_N$.\\

\noindent\textbf{Remark:} Note that the increment $\Delta X_{H,N}(s)$ is a weighted sum of increments
of the Bernoulli random walk $B_N$ --- it is a linear combination of the current coin toss $\omega_N(s)$ 
and all previous coin tosses $\omega_N(r)$, $r<s$. Observe also that since 
$\lim_{H\uparrow\frac{1}{2}}-(H-\frac{1}{2})\zeta(\frac{3}{2}-H)=\lim_{s\downarrow 1}(s-1)\zeta(s)=1$,
the two cases meet continuously at $H=\frac{1}{2}$. For $H>\frac{1}{2}$, we may actually choose $K_{H}$
as we please since the term will vanish in the limit (see below), but $K_{H}=1$ is the natural value and probably the one that
gives best results in numerical work. \\

We are now ready to state the main result. Note that when $H=\frac{1}{2}$, $\Delta X_{\frac{1}{2},N}(t)=\Delta B_N(t)$ and the theorem just reduces to the classical convergence of a Bernoulli random walk to Brownian motion.

\begin{theorem}[Main Theorem] For all real numbers $H$, $0<H<1$, the processes $c_{H}X_{H,N}$ converge weakly in $D([0,\infty))$ to fractional
Brownian motion with Hurst index $H$. 
\end{theorem}

\noindent\textbf{Notation:} In the rest of the paper, we drop the notational
dependence on $N$ and $H$, and write simply $X$, $B$, $T$, $\Delta t$ for $X_{H,N}$, $B_N$, $T_N$, $\Delta t_N$ etc. when no confusion can arise.\\

As we are interested in understanding the dynamics of fractional Brownian motion, we have defined $X$ by specifying its increments $\Delta X(s)$. To prove
the main theorem, we need an expression for $X(t)$. This is just a small calculation:

\[X(t)=\sum_{s=0}^t \Delta X(s)=\sum_{s=0}^t K_H\Delta t^{H-\frac{1}{2}}\Delta B_s
+\sum_{s=0}^t \sum_{r=-\infty}^s (H-\frac{1}{2})(s-r)^{H-\frac{3}{2}}\Delta t\Delta B_r\]
Changing the order of summation, we have
\[X(t)=K_H\Delta t^{H-\frac{1}{2}} B_t
+\sum_{r=0}^t\sum_{s=r+\Delta t}^t (H-\frac{1}{2})(s-r)^{H-\frac{3}{2}}\Delta t \Delta B_r\]
\[+\sum_{r=-\infty}^0 \sum_{s=0}^t (H-\frac{1}{2})(s-r)^{H-\frac{3}{2}}\Delta t\Delta B_r\]
where $B_t=\sum_{r=0}^t\Delta B_r$ is a random walk converging to Brownian motion. Observe that when $H>\frac{1}{2}$, the first term $K_H\Delta t^{H-\frac{1}{2}} B_t$ vanishes when $N\to\infty$ (this is why the choice of $K_H$ is irrelevant in this case), but when $H<\frac{1}{2}$, the term explodes. In this case we have a delicate balance between two terms going to infinity, and a correct choice of $K_H$ is crucial.

The idea is now to simplify the expression for $X$ by replacing the sums $\sum(H-\frac{1}{2})(s-r)^{H-\frac{3}{2}}\Delta t$ by the corresponding integrals $\int (H-\frac{1}{2})(s-r)^{H-\frac{3}{2}}\; ds$, and then performing the integration. This works
nicely for $H>\frac{1}{2}$, but when $H<\frac{1}{2}$, one of the integrals diverges, and we have to be more careful. Put crudely, it is the divergence of this integral that will cancel the divergence of the term $K_H\Delta t^{H-\frac{1}{2}} B_t$.

We are ready to prove the main theorem, and start with the simplest case.

\section{The case $H>\frac{1}{2}$}
We start from the expression 
\[X(t)=\Delta t^{H-\frac{1}{2}} B_t
+\sum_{r=0}^t\sum_{s=r+\Delta t}^t (H-\frac{1}{2})(s-r)^{H-\frac{3}{2}}\Delta t \Delta B_r+\]
\[+\sum_{r=-\infty}^0 \sum_{s=0}^t (H-\frac{1}{2})(s-r)^{H-\frac{3}{2}}\Delta t\Delta B_r\]
above (remember that $K_{H}=1$ in this case). Since $H>\frac{1}{2}$, we have no problem
with convergence, and if we let $\epsilon_N(r,t)$ be the error term:
\[\epsilon_N(r,t):=\sum_{s=r+\Delta t}^t (H-\frac{1}{2})(s-r)^{H-\frac{3}{2}}\Delta t -
\int_{r+\Delta t}^t  (H-\frac{1}{2})(s-r)^{H-\frac{3}{2}}\;ds,\]
we get
\[\sum_{s=r+\Delta t}^t (H-\frac{1}{2})(s-r)^{H-\frac{3}{2}}\Delta t =
\int_{r+\Delta t}^t  (H-\frac{1}{2})(s-r)^{H-\frac{3}{2}}\;ds+\epsilon_{N}(r,t)=\]
\[=(t-r)^{H-\frac{1}{2}}-\Delta t^{H-\frac{1}{2}}+\epsilon_{N}(r,t)\]
Similarly, with 
\[\delta_N(r,t):=\sum_{s=0}^t (H-\frac{1}{2})(s-r)^{H-\frac{3}{2}}\Delta t -
\int_{0}^t  (H-\frac{1}{2})(s-r)^{H-\frac{3}{2}}\;ds,\]
we get
\[\sum_{s=0}^t (H-\frac{1}{2})(s-r)^{H-\frac{3}{2}}\Delta t =
\int_{0}^t  (H-\frac{1}{2})(s-r)^{H-\frac{3}{2}}\;ds+\delta_{N}(r,t)=\]
\[=(t-r)^{H-\frac{1}{2}}-(-r)^{H-\frac{1}{2}}+\delta_{N}(r,t)\]
This means that

\[X(t)=\sum_{r=0}^t\left((t-r)^{H-\frac{1}{2}}+\epsilon_{N}(r)\right) \Delta B_r+\]
\[+\sum_{r=-\infty}^0 \left((t-r)^{H-\frac{1}{2}}-(-r)^{H-\frac{1}{2}}+\delta_{N}(r)\right)\Delta B_r=\]
\[=\sum_{r=-\infty}^t\left( (t-r)^{H-\frac{1}{2}}- (-r)_{+}^{H-\frac{1}{2}}\right)\Delta B_{r}
+\sum_{r=0}^t\epsilon_{N}(r,t) \Delta B_r+
\sum_{r=-\infty}^0  \delta_{N}(r,t)\Delta B_r\]
We want to prove that $X$ converges weakly to fractional Brownian motion. According to Theorem 1 in \cite{KS}, it suffices to show that $E(c_H^2X(t)^2)\to t^{2H}$. This follows immediately from the Mandelbrot-Van Ness representation and the following lemma.

 \begin{lemma} For $\frac{1}{2}<H<1$:
 \begin{itemize}
 \item[(i)] $E\left((\sum_{r=0}^t\epsilon_N(r,t)\Delta B_r)^2\right)\leq (H-\frac{1}{2})^2t\Delta t^{2H-1}$
  \item[(ii)] $E\left((\sum_{r=-\infty}^0\delta_N(r,t)\Delta B_r)^2\right)\leq (H-\frac{1}{2})^2\zeta(3-2H)\Delta t^{2H}$
 \end{itemize}
 \end{lemma}
 \noindent\emph{Proof:} (i) We first observe that
 \[\epsilon_N(r,t)=\sum_{s=r+\Delta t}^t (H-\frac{1}{2})(s-r)^{H-\frac{3}{2}}\Delta t -
\int_{r+\Delta t}^t  (H-\frac{1}{2})(s-r)^{H-\frac{3}{2}}\;ds>0\]
since $\sum_{s=r+\Delta t}^t (H-\frac{1}{2})(s-r)^{H-\frac{3}{2}}\Delta t$ is an upper Riemann sum for the integral. Since $\sum_{s=r+2\Delta t}^t (H-\frac{1}{2})(s-r)^{H-\frac{3}{2}}\Delta t$ is a lower Riemann sum, we also have
\[0\leq\epsilon_N(r,t)\leq (H-\frac{1}{2})\Delta t^{H-\frac{3}{2}}\Delta t=(H-\frac{1}{2})\Delta t^{H-\frac{1}{2}}\]
Thus
\[E\left((\sum_{r=0}^t\epsilon_N(r,t)\Delta B_r)^2\right)=
\sum_{r=0}^t\epsilon_N(r,t)^2\Delta t\leq\]
\[\leq \sum_{r=0}^t(H-\frac{1}{2})^2\Delta t^{2H-1}\Delta t\leq (H-\frac{1}{2})^2t\Delta t^{2H-1}\]

(ii) Using approximating Riemann sums as in part (i), we see that
\[0\leq\delta_N(r,t)\leq (H-\frac{1}{2})(-r)^{H-\frac{3}{2}}\Delta t,\]
and thus
\[E\left((\sum_{r=-\infty}^0\delta_N(r,t)\Delta B_r)^2\right)=
\sum_{-\infty}^0\delta_N(r,t)^2\Delta t\leq \sum_{r=-\infty}^0(H-\frac{1}{2})^2(-r)^{2H-3}\Delta t^3\]
Letting $r=-k\Delta t$, we get
\[E\left((\sum_{r=-\infty}^0\delta_N(r,t)\Delta B_r)^2\right)\leq \sum_{k=0}^{\infty}(H-\frac{1}{2})^2k^{2H-3}\Delta t^{2H}=\]
\[=(H-\frac{1}{2})^2\zeta(3-2H)\Delta t^{2H}\]
This completes the proof of the lemma (and also the proof of the Main Theorem for the case $H>\frac{1}{2}$). \QED

\section{The case $H<\frac{1}{2}$}
Again we start from the expression 
\[X(t)=K_{H}\Delta t^{H-\frac{1}{2}} B_t
+\sum_{r=0}^t\sum_{s=r+\Delta t}^t (H-\frac{1}{2})(s-r)^{H-\frac{3}{2}}\Delta t \Delta B_r+\]
\[+\sum_{r=-\infty}^0 \sum_{s=0}^t (H-\frac{1}{2})(s-r)^{H-\frac{3}{2}}\Delta t\Delta B_r\]
In this case, one of the integrals we worked with above diverges, and we have to be more careful.
Let us start with a closer look at the term $\sum_{s=r+\Delta t}^t (H-\frac{1}{2})(s-r)^{H-\frac{3}{2}}\Delta t$. We obviously have
\[\sum_{s=r+\Delta t}^t (H-\frac{1}{2})(s-r)^{H-\frac{3}{2}}\Delta t=\]
\[\sum_{s=r+\Delta t}^{\infty} (H-\frac{1}{2})(s-r)^{H-\frac{3}{2}}\Delta t
-\sum_{s=t}^{\infty} (H-\frac{1}{2})(s-r)^{H-\frac{3}{2}}\Delta t\]
and if we let $r=N\Delta t$, $s=k\Delta t$, we get
\[\sum_{s=r+\Delta t}^{\infty} (H-\frac{1}{2})(s-r)^{H-\frac{3}{2}}\Delta t=\sum_{k=N+1}^{\infty}(H-\frac{1}{2})(k\Delta t-N\Delta t)^{H-\frac{3}{2}}\Delta t\]
\[=(H-\frac{1}{2})\Delta t^{H-\frac{1}{2}}\sum_{k=N+1}^{\infty}(k-N)^{H-\frac{3}{2}}=
(H-\frac{1}{2})\Delta t^{H-\frac{1}{2}}\sum_{n=1}^{\infty}n^{H-\frac{3}{2}}\]
\[=(H-\frac{1}{2})\Delta t^{H-\frac{1}{2}}\zeta(\frac{3}{2}-H)=-K_{H}\Delta t^{H-\frac{1}{2}}\]
Substituting this into the expression for $X(t)$, we get
\[X(t)=\sum_{r=0}^t\sum_{s=t}^{\infty} -(H-\frac{1}{2})(s-r)^{H-\frac{3}{2}}\Delta t\Delta B_r\]
\[+\sum_{r=\infty}^0 \sum_{s=0}^t (H-\frac{1}{2})(s-r)^{H-\frac{3}{2}}\Delta t\Delta B_r\]
The two sums in this expression have less dangerous limits than the one we just got rid of, and
can be approximated by integrals. If we let
\[\tilde{\epsilon}_N(r,t):=\sum_{s=t}^{\infty} -(H-\frac{1}{2})(s-r)^{H-\frac{3}{2}}\Delta t-\int_t^{\infty}-(H-\frac{1}{2})(s-r)^{H-\frac{3}{2}}\;ds,\]
we get  (remember
that $H<\frac{1}{2}$):
\[\sum_{s=t}^{\infty} -(H-\frac{1}{2})(s-r)^{H-\frac{3}{2}}\Delta t=
\int_t^{\infty}-(H-\frac{1}{2})(s-r)^{H-\frac{3}{2}}\;ds+\tilde{\epsilon}_N(r,t)\]
\[=\left[-(s-r)^{H-\frac{1}{2}}\right]_{s=t}^{s=\infty}+\tilde{\epsilon}_N(r,t)
=(t-r)^{H-\frac{1}{2}}+\tilde{\epsilon}_N(r,t)\]
Similarly, if we let
\[\tilde{\delta}_N(r,t):=\sum_{s=0}^{t}- (H-\frac{1}{2})(s-r)^{H-\frac{3}{2}}\Delta t-
\int_0^{t}-(H-\frac{1}{2})(s-r)^{H-\frac{3}{2}}\;ds,\]
we get
\[\sum_{s=0}^{t} (H-\frac{1}{2})(s-r)^{H-\frac{3}{2}}\Delta t=
\int_0^{t}(H-\frac{1}{2})(s-r)^{H-\frac{3}{2}}\;ds-\tilde{\delta}_N(r,t)\]
\[=\left[(s-r)^{H-\frac{1}{2}}\right]_{s=0}^{s=t}-\tilde{\delta}_N(r,t)
=(t-r)^{H-\frac{1}{2}}-(-r)^{H-\frac{1}{2}}-\tilde{\delta}_N(r,t)\]
We thus have
\[X(t)= \sum_{r=0}^t\left((t-r)^{H-\frac{1}{2}}+\tilde{\epsilon}_N(r,t)\right)\Delta B_r+\]
\[+\sum_{r=-\infty}^0\left((t-r)^{H-\frac{1}{2}}-(-r)^{H-\frac{1}{2}}-\tilde{\delta}_N(r,t)\right) \Delta B_r\]
\[= \sum_{r=-\infty}^t \left((t-r)^{H-\frac{1}{2}}-(-r)_+^{H-\frac{1}{2}} \right)\Delta B_r+\]
\[+\sum_{r=0}^t\tilde{\epsilon}_N(r,t)\Delta B_r
-\sum_{r=-\infty}^{0}\tilde{\delta}_N(r,t)\Delta B_r\]
To prove that $c_HX$ converges weakly to fractional Brownian motion, we can now longer use Theorem 1 of \cite{KS} as in the previous case since this theorem requires that $H>\frac{1}{2}$. However, the first term in the expression above obviously converges weakly to 
\[\int_{r=-\infty}^t \left((t-r)^{H-\frac{1}{2}}-(-r)_+^{H-\frac{1}{2}} \right)\;db_r,\]
 and the next lemma shows that error terms go uniformly to zero. Using the Mandelbrot-Van Ness representation, we then get the Main Theorem for $H<\frac{1}{2}$.
 
 \begin{lemma}For each $H$, $0<H<\frac{1}{2}$, there is a constant $K_H\in\real_+$ (independent of $N$ and $t$) such that
  \[\left|X(t)- \sum_{r=-\infty}^t \left((t-r)^{H-\frac{1}{2}}-(-r)_+^{H-\frac{1}{2}}\right) \Delta B_r\right| \leq K_H\Delta t^{H}\]
  \end{lemma}
  \noindent\emph{Proof:} It clearly suffices to show that there are constants $C_H, D_H\in\real_+$ (independent of $N$ and $t$) such that
  \[\left|\sum_{r=0}^t\tilde{\epsilon}_N(r,t)\Delta B_r\right|\leq C_H\Delta t^{H}
  \hspace{5mm}\mbox{and}\hspace{5mm}\left|\sum_{r=-\infty}^{0}\tilde{\delta}_N(r)\Delta B_r\right|\leq D_H\Delta t^{H}\]

We begin with the $\tilde{\epsilon}_N$-case. By definition
  \[\tilde{\epsilon}_N(r,t)=\sum_{s=t}^{\infty} -(H-\frac{1}{2})(s-r)^{H-\frac{3}{2}}\Delta t-
\int_t^{\infty}-(H-\frac{1}{2})(s-r)^{H-\frac{3}{2}}\;ds\]
Since $\sum_{s=t}^{\infty} -(H-\frac{1}{2})(s-r)^{H-\frac{3}{2}}\Delta t$ is an upper
Riemann sum for the integral
$\int_t^{\infty}-(H-\frac{1}{2})(s-r)^{H-\frac{3}{2}}\;ds$, and $\sum_{s=t+\Delta t}^{\infty} -(H-\frac{1}{2})(s-r)^{H-\frac{3}{2}}\Delta t$ is a lower Riemann sum, we have 
\[0\leq \tilde{\epsilon}_N(r,t)\leq-(H-\frac{1}{2})(t-r)^{H-\frac{3}{2}}\Delta t\]
Hence (remember that $|\Delta B_r|=\Delta t^{\frac{1}{2}}$)
  \[\left|\sum_{r=0}^t\tilde{\epsilon}_N(r,t)\Delta B_r\right|
  \leq \sum_{r=0}^t-(H-\frac{1}{2})(t-r)^{H-\frac{3}{2}}\Delta t^{\frac{3}{2}}\]
  If we let $t=K\Delta t$, $r=k\Delta t$, we can rewrite the last sum as
 \[\sum_{k=0}^{K-1}-(H-\frac{1}{2})(K-k)^{H-\frac{3}{2}}\Delta t^{H} \leq -(H-\frac{1}{2})\zeta(\frac{3}{2}-H)\Delta t^{H}\]
 This completes the $\tilde{\epsilon}_N$-part of the argument.
 
 Turning to the term $\sum_{r=-\infty}^{0}\tilde{\delta}_N(r)\Delta B_r$, we first observe that by definition
 \[\tilde{\delta}_N(r)=\sum_{s=0}^{t} -(H-\frac{1}{2})(s-r)^{H-\frac{3}{2}}\Delta t-
\int_0^{t}-(H-\frac{1}{2})(s-r)^{H-\frac{3}{2}}\;ds\]
Again, $\sum_{s=0}^{t} -(H-\frac{1}{2})(s-r)^{H-\frac{3}{2}}\Delta t$ is an upper
Riemann sum, and we easily see that
\[0\leq \tilde{\delta}_N(r)\leq  -(H-\frac{1}{2})(-r)^{H-\frac{3}{2}}\Delta t\]
Letting $r=-k\Delta t$, we get
 \[E\left|\sum_{r=-\infty}^0\tilde{\delta}_N(r)\Delta B_r\right|
   \leq \sum_{r=-\infty}^0-(H-\frac{1}{2})(-r)^{H-\frac{3}{2}}\Delta t^{\frac{3}{2}}\leq\]
  \[\leq -(H-\frac{1}{2})\Delta t^{H}\sum_{k=0}^{\infty}k^{H-\frac{3}{2}}=-(H-\frac{1}{2})\zeta(\frac{3}{2}-H)\Delta t^{H}\]
  This proves the lemma (and hence the Main Theorem for the remaining case $H<\frac{1}{2}$). \QED

\end{document}